\documentclass[review]{elsarticle}

\usepackage{amsmath,amssymb,amsxtra,comment,graphicx,psfrag}
\usepackage{bm,mathrsfs}
\usepackage{mathtools}
\usepackage{xspace}
\usepackage{color}
\usepackage{empheq}
\usepackage{enumerate,enumitem}
\usepackage{hhline}
\usepackage{fancyhdr}
\usepackage{pgfplots}

\usepackage{pifont}
\usepackage{textcomp}
\usepackage[english]{babel}
\usepackage{framed}
\usepackage{bm}
\usepackage{subfigure}
\usepackage{array}
\usepackage{booktabs}
\usepackage{cases}
\biboptions{numbers, sort&compress}

\newtheorem{theorem}{Theorem}

\journal{Journal of Computational Physics}

\def\x{{\bf x}}

\def\d{{\mathrm d}}

\def\R{{\mathbb R}}

\begin{document}

\begin{frontmatter}

\title{A perimeter-decreasing and area-conserving algorithm for surface diffusion flow of curves}
\tnotetext[mytitlenote]{This work is supported by the National Natural Science Foundation of China (project no. 11871384) and the Hong Kong Research Grants Council (GRF project no.: 15300920).}


\author[mymainaddress]{Wei Jiang}
\ead{jiangwei1007@whu.edu.cn}

\author[mysecondaryaddress]{Buyang Li\corref{mycorrespondingauthor}}
\cortext[mycorrespondingauthor]{Corresponding author}
\ead{buyang.li@polyu.edu.hk}

\address[mymainaddress]{School of Mathematics and Statistics $\&$ Hubei Key Laboratory of Computational Science, Wuhan University, Wuhan 430072, P. R. China.}
\address[mysecondaryaddress]{Department of Applied Mathematics, The Hong Kong Polytechnic University, Hong Kong.}

\begin{abstract}
A fully discrete finite element method, based on a new weak formulation and a new time-stepping scheme, is proposed for the surface diffusion flow of closed curves in the two-dimensional plane. It is proved that the proposed method can preserve two geometric structures simultaneously at the discrete level, i.e.,  the perimeter of the curve decreases in time while the area enclosed by the curve is conserved. Numerical examples are provided to demonstrate the convergence of the proposed method and the effectiveness of the method in preserving the two geometric structures.
\end{abstract}


\end{frontmatter}



\section{Introduction}\label{sec:intro}

This article concerns the numerical approximation to the surface diffusion flow of closed curves in the two-dimensional plane, i.e., the evolution of a curve
$$
\Gamma[\mathbf X(\cdot,t)] = \{\mathbf X(\xi,t): \xi \in I \}, \quad t\in[0,T],
$$
determined by a parametrization $\mathbf X(\cdot,t):I\rightarrow \R^2$, where $I=[0,1]$ is the periodic unit interval  (one-dimensional torus which identifies $0$ and $1$), satisfying the fourth-order geometric evolution equation
\begin{equation}\label{SF_Eq}
\left\{
\begin{aligned}
     &\partial_t\mathbf{X} = - [\partial_{\Gamma}^2 (\boldsymbol\kappa\cdot \boldsymbol\nu)]\boldsymbol\nu,\\
     &\boldsymbol\kappa = \partial_{\Gamma} \boldsymbol{\tau},
\end{aligned}
\right.  \qquad \text{on} \quad   I \times [0,T],
\end{equation}
where $\boldsymbol\kappa=\boldsymbol\kappa(\xi, t)$, $\boldsymbol\nu= \boldsymbol\nu(\xi, t)$ and $\boldsymbol\tau= \boldsymbol\tau(\xi, t)$ are the curvature vector, inward unit normal vector and unit tangential vector of the curve $\Gamma[\mathbf X(\cdot,t)]$ at the point $\mathbf X(\xi,t)$. The notation $\partial_\Gamma$ represents tangential differentiation along the curve $\Gamma$, defined by
$$
\partial_\Gamma v(\xi,t) := \frac{\partial_\xi v(\xi,t)}{|\partial_\xi \mathbf X(\xi,t)|} \quad\mbox{for a function $v$ defined on $I\times[0,T]$} .
$$

The parametric equation \eqref{SF_Eq} for surface diffusion flow was proposed by Mullins in 1957 for modeling the evolution of microstructure in polycrystalline materials~\cite{Mullins57}. It was generalized to include anisotropic effects of crystalline films in  \cite{Cahn94,Gurtin02,LiBo09,Torabi10}. These models play important roles in various applications in materials science and solid-state dewetting; see \cite{Thompson12,Jiang12,Wang15,Jiang16,Bao17b,Jiang18a,Jiang18b,Jiang19c} and the references therein.

It is well-known that the surface diffusion flow described by \eqref{SF_Eq} is the $H^{-1}$ gradient flow of the perimeter functional (see \cite{Cahn94,Taylor94})
$$
|\Gamma[\mathbf X(\cdot,t)]|
=\int_I |\partial_\xi \mathbf X(\xi,t)|\d \xi  .
$$
As a result, the surface diffusion flow has two geometric structures:
\begin{itemize}
\item[{(i)}]The perimeter of the evolving curve decreases in time;
\item[{(ii)}]The area enclosed by the curve is conserved in the evolution.
\end{itemize}
These geometric structures are also desired in the numerical approximation of surface diffusion flow.

Many numerical methods have been developed for simulating the evolution of a curve/surface in 2D/3D governed by surface diffusion flow and related geometric evolution equations. 
For the non-parametric equations describing axially symmetric surface diffusion flow, Coleman, Falk \& Moakher \cite{Coleman-Falk-Moakher-1996} proposed space-time finite element methods (FEMs) to preserve the perimeter decrease and area conservation properties, and Deckelnick, Dziuk \& Elliott \cite{Deckelnick-Dziuk-Elliott-2003} proved convergence of semidiscrete FEMs. 
For the non-parametric equations describing surface diffusion flow of graphs, variational formulation, error analysis and numerical simulation were done by B\"ansch, Morin \& Nochetto \cite{Bansch04} and 
Deckelnick, Dziuk \& Elliott \cite{Deckelnick-Dziuk-Elliott-2005-SINUM} for semidiscrete and fully discrete FEMs, respectively. 
In the axially symmetric and graph cases, both semi-discrete and fully discrete methods can preserve the perimeter decrease and area conservation. 

For general parametric surface diffusion flow, the underlying numerical methods are the parametric FEMs developed by Dziuk \cite{Dziuk90,Dziuk1999}, Deckelnick, Dziuk \& Elliott \cite{Deckelnick-Dziuk-Elliott-2005}, B\"ansch, Morin \& Nochetto \cite{Bansch05}, Dziuk \& Elliott \cite{Dziuk-Elliott-2007}, Barrett, Garcke \& N\"{u}rnberg \cite{Barrett07a,Barrett08b}, etc. 
In particular, B\"ansch, Morin \& Nochetto \cite{Bansch05} proposed the first parametric FEM for the surface diffusion flow. The method preserves the decrease of perimeter (see \cite[Theorem 2.1]{Bansch05}), while the area enclosed by the curve may change in time. Parametric FEMs based on novel variational formulations were developed in \cite{Barrett07a,Barrett08b,Barrett11,Barrett19b}. These methods yield good mesh distribution and unconditional stability. The generalization to anisotropic surface diffusion flows and axisymmetric geometry cases was made in \cite{Barrett08,Barrett08a,Barrett19}. 
These methods can preserve the decrease of perimeter in 2D (or decrease of area in 3D) in the fully-discrete finite element scheme, and the area conservation in 2D (or volume conservation in 3D) only in the semi-discrete FEM.
More recently, Jiang et al. developed parametric FEMs for simulating solid-state dewetting problems described by surface diffusion flow and contact line migration in 2D~\cite{Bao17} and 3D~\cite{Jiang19b}, without considering the perimeter/area decrease and area/volume conservation properties in the discrete level.

Overall, existing numerical methods for parametric surface diffusion flow often yield curves/surfaces with decreasing perimeter/area, while the area/volume enclosed by the curve/surface is conserved only in the semi-discrete FEMs. In this paper, we construct a fully discrete parametric FEM for surface diffusion flow of closed curves, based on a new weak formulation and a new time-stepping scheme, to preserve the two geometric structures simultaneously.

The rest of paper is organized as follows. In section 2, we propose a new weak formulation for solving  problem \eqref{SF_Eq}, then a parametric FEM is used to discretize the weak formulation and results in a nonlinear system of equations. In section 3, we rigorously prove that the numerical scheme is area-preserving and perimeter-decreaseing. In section 4, we present some numerical simulations to demonstrate the accuracy and efficiency of the proposed method. Finally, we draw some conclusions in section 5.



\section{The numerical method}

In this section, we introduce a new variational formulation of the surface diffusion flow equation \eqref{SF_Eq} and then propose fully discrete FEMs corresponding to the variational forms.

\subsection{Notation of spatial and temporal discretizations}
Let $0 = \xi_0<\xi_1 \cdots < \xi_M = 1$ be a quasi-uniform partition of the interval $I$ (where we identify $\xi_0$ and $\xi_M$), which is divided into $M$ subintervals labeled as $I_j := [\xi_{j - 1}, \xi_j]$, $j = 1,\cdots, M$. We define the following finite element spaces of piecewise linear functions subject to the partition, i.e.,
\begin{align}
    & V_h = \big\{\varphi\in C(I; \mathbb{R}): \varphi(0) = \varphi(1), ~\varphi|_{I_j}\in P_1, ~j = 1, 2, \cdots, M \big\},\\
    & {\bf V}_h = \big\{\boldsymbol\psi\in C(I; \mathbb{R}^2): \boldsymbol\psi(0) = \boldsymbol\psi(1), ~\boldsymbol\psi|_{I_j}\in P_1 \times P_1, ~j = 1, 2, \cdots, M \big\},
\end{align}
where $P_1$ denotes the space of polynomials of degree $\le 1$.

Let $t_n=n\tau$, $n=0,1,2,\dots,N$, be a uniform partition of the time interval $[0,T]$ with stepsize $\tau=T/N$. For each $n=1,\dots,N$, we denote by $\Gamma_h^{n}$ a piecewise linear curve determined by a sequence of non-coincident nodes $\mathbf{x}_j^{n}\in \mathbb{R}^2$, $j = 1, \cdots, M$ in the counter clockwise order, approximating the closed curve $\Gamma[{\bf X}(\cdot,t_n)]$. Correspondingly, the nodal vector $\mathbf{x}^n = (\mathbf{x}_1^n, \cdots, \mathbf{x}_M^n)$ uniquely determines a finite element function $\mathbf{X}_h^n\in {\bf V}_h$ which parametrizes $\Gamma_h^n$, satisfying
$$
\mathbf{X}_h^n(\xi_j) = \mathbf{x}_j^n,\quad j = 1, \cdots, M.
$$
We assume that $\mathbf{X}_h^{n-1}\in {\bf V}_h$ is given and look for $\mathbf{X}_h^{n}\in {\bf V}_h$ by an implicit time-stepping scheme defined below.

In the time interval $[t_{n-1}, t_{n}]$, we define the intermediate curves
$$\Gamma_h(t) = \{\mathbf{X}_h(\xi, t): \xi\in[0, 1]\}$$ with parametrization
\begin{equation}\label{eq:u_intermediate}
    \mathbf{X}_h(\xi, t) = \dfrac{t_{n} - t}{\tau}\mathbf{X}_h^{n - 1}(\xi) + \dfrac{t - t_{n-1}}{\tau}\mathbf{X}_h^n(\xi), \qquad t\in [t_{n-1}, t_{n}] .
\end{equation}
Then for a fixed $\xi\in I$ the material point $\mathbf{X}(\xi, t)$ on the intermediate curve $\Gamma_h(t)$ moves with velocity
\begin{equation}\label{eq:velocity}
    \mathbf{v}_h^n(\xi) = \dfrac{\mathbf{X}_h^n(\xi) - \mathbf{X}_h^{n - 1}(\xi)}{\tau},
\end{equation}
which is a finite element function in ${\bf V}_h$ independent of time.

The unit normal and tangential vectors at the $\mathbf{X}_h(\xi, t)\in \Gamma_h(t)$ are given by
\begin{equation}\label{eq:nu_tau_intermediate}
    \boldsymbol\nu_h(\xi, t) = \dfrac{\partial_\xi\mathbf{X}_h(\xi, t)^\perp}{|\partial_\xi\mathbf{X}_h(\xi, t)|}
    \quad\mbox{and} \quad
    \boldsymbol\tau_h(\xi, t) = \dfrac{\partial_\xi\mathbf{X}_h(\xi, t)}{|\partial_\xi\mathbf{X}_h(\xi, t)|} ,
\end{equation}
respectively, where $^\perp$ denotes rotation of a vector by an angle $\pi/2$ in the counter clockwise direction. These notations will be used in the following two subsections to define numerical schemes for \eqref{SF_Eq}.



\subsection{A structure preserving FEM}\label{section:I}

We decompose the curvature vector into its normal and tangential components, separately, i.e.,
$$
\boldsymbol\kappa= p\,\boldsymbol\nu + q\,\boldsymbol\tau,
$$
where $p$ and $q$ are scalar functions (in fact, $q=0$ for the continuous problem).
With these notations, we rewrite equation \eqref{SF_Eq} as
\begin{subequations}\label{eq:surface_diffusion_rewrite}
\begin{align}
\partial_t\mathbf{X} \cdot \boldsymbol \nu &= - \partial_{\Gamma}^2 p,\\
  \partial_t\mathbf{X} \cdot \boldsymbol \tau &= 0,\\
  p\,\boldsymbol\nu + q\,\boldsymbol\tau &= \partial_{\Gamma} \boldsymbol\tau .
\end{align}
\end{subequations}

Let $H^1(I)=\{\phi \in H^1(0,1): \phi(0) = \phi(1)\}$ and consider the following weak formulation of \eqref{eq:surface_diffusion_rewrite}:
find $\mathbf{X}(\cdot, t)\in H^1(I)\times H^1(I)$ and $p(\cdot, t),q(\cdot, t)\in H^1(I)$ such that the equations
\begin{subequations}\label{eq:variation_formulation_conti}
\begin{align}
&\hspace{-10pt} \int_{I} \partial_t\mathbf{X}(\cdot,t) \cdot (\phi\boldsymbol\nu) \, \mathrm{d}\Gamma(t) = \int_{I} \partial_\Gamma p\, \partial_\Gamma \phi \,\mathrm{d}\Gamma(t), &&\forall\phi\in H^1(I), \\
&\hspace{-10pt} \int_{I} \partial_t \mathbf{X}(\cdot,t) \cdot (\varphi\boldsymbol\tau) \, \mathrm{d}\Gamma(t) = 0, && \forall\varphi\in H^1(I),\\
&\hspace{-10pt} \int_{I} (p\,\boldsymbol\nu + q\,\boldsymbol\tau) \cdot \boldsymbol\psi \,\mathrm{d}\Gamma(t) = \int_{I} \partial_\Gamma \boldsymbol{\tau} \cdot \boldsymbol\psi \,\mathrm{d} \Gamma(t), && \forall\boldsymbol\psi\in H^1(I)\times H^1(I) ,
\end{align}
\end{subequations}
hold for all $t\in(0,T]$ under an initial condition $\mathbf{X}(\xi, 0)=\mathbf{X}_0(\xi)$, where $\mathbf{X}_0$ is a given parametrization of the curve at time $t=0$, and
$$\d\Gamma(t):=|\partial_\xi {\bf X}(\xi,t)|\d\xi .$$

On the intermediate curve $\Gamma_h(t)$, we approximate the curvature at the point $\mathbf{X}_h(\xi, t)$ by a function
\begin{equation}\label{eq:kappa}
    \boldsymbol\kappa_h(\xi, t) = p_h(\xi)\,\boldsymbol\nu_h(\xi, t) + q_h(\xi)\,\boldsymbol\tau_h(\xi, t),
\end{equation}
with some finite element functions $p_h, q_h\in V_h$ independent of time when $t\in [t_{n-1},t_n]$.


By using the discrete velocity and curvature defined by \eqref{eq:velocity} and \eqref{eq:kappa}, we look for $\mathbf{X}_h^n\in {\bf V}_h$ and $p_h^n, q_h^n\in V_h$ satisfying the following equations:
\begin{subequations}\label{eq:variation_formulation_initial}
\begin{align}
&\hspace{-10pt}  \int_{t_{n-1}}^{t_{n}}\int_{I} \mathbf{v}_h^n\cdot( \phi_h\, \boldsymbol\nu_h)\, \d\Gamma_h(t) \mathrm{d}t      =\int_{t_{n-1}}^{t_{n}}\int_{I}
\partial_\Gamma p_h^n\, \partial_\Gamma \phi_h\, \d\Gamma_h(t) \mathrm{d}t, && \forall \phi_h \in V_h, \label{eq:variation_formulation_initial_a}\\
&\hspace{-10pt}\int_{t_{n-1}}^{t_{n}}\int_{I} \mathbf{v}_h^n\cdot ( \varphi_h\,\boldsymbol\tau_h) \, \d\Gamma_h(t)\mathrm{d}t
= 0,
&& \forall \varphi_h\in V_h, \label{eq:variation_formulation_initial_b}\\
&\hspace{-10pt}     \int_{t_{n-1}}^{t_{n}} \hspace{-2pt}  \int_{I}  ( p_h^n\,\boldsymbol\nu_h +  q_h^n\, \boldsymbol\tau_h)\cdot \boldsymbol\psi_h \,\d\Gamma_h(t)\mathrm{d}t
= \sum_{j = 1}^M\boldsymbol\psi_h(\xi_j) \cdot  \int_{t_{n-1}}^{t_{n}}[\boldsymbol\tau_h]_j \,\mathrm{d}t,  && \forall\boldsymbol\psi_h\in {\bf V}_h, \label{eq:variation_formulation_initial_c}
\end{align}
\end{subequations}
where $\d\Gamma_h(t):=|\partial_\xi {\bf X}_h(\xi,t)|\d\xi$ and
$$
[\boldsymbol\tau_h]_j = \lim_{\varepsilon\to 0}\left(\boldsymbol\tau_h(\xi_j + \varepsilon, t) - \boldsymbol\tau_h(\xi_j - \varepsilon, t)\right)
$$ denotes the jump of tangential vector at the node $\mathbf{X}_h(\xi_j ,t)$. The right side of \eqref{eq:variation_formulation_initial_c} is obtained by using the identity
$$
\partial_\xi \boldsymbol \tau_h(\xi,t)
= \sum_{j=1}^M \, [\boldsymbol\tau_h]_j \, \delta(\xi-\xi_j)
$$
for the piecewise constant function $\boldsymbol \tau_h$,
where $\delta(\xi-\xi_j)$ denotes the Dirac delta function centered at the point $\xi_j$.

If we denote
\begin{align}
    & \boldsymbol\nu_h^n(\xi) = \int_{t_{n-1}}^{t_{n}} \boldsymbol\nu_h(\xi, t) |\partial_{\xi} \mathbf{X}_h(\xi, t) | \, \mathrm{d}t = \int_{t_{n-1}}^{t_{n}} (\partial_{\xi}\mathbf{X}_h(\xi, t) )^{\perp}\, \mathrm{d}t, \label{eq:nu_tau_int} \\
    & \boldsymbol\tau_h^n(\xi)  = \int_{t_{ n -1}}^{t_{n}}\boldsymbol\tau_h(\xi, t) |\partial_{\xi}\mathbf{X}_h(\xi, t) |\, \mathrm{d}t = \int_{t_{ n - 1}}^{t_{n}}\partial_{\xi} \mathbf{X}_h(\xi, t) \, \mathrm{d}t, \label{eq:tau_int}\\
    & \rho_h^n(\xi) = \int_{t_{n-1}}^{t_{n}} |\partial_{\xi}\mathbf{X}_h(\xi, t) |^{-1}\, \mathrm{d}t, \label{eq:rho_int}
\end{align}
then \eqref{eq:variation_formulation_initial} can be equivalently as:
find $\mathbf{X}_h^n\in {\bf V}_h$ and $p_h^n, q_h^n \in V_h$ such that
\begin{subequations}\label{eq:variation_formulation}
\begin{align}
& \int_I \dfrac{\mathbf{X}_h^n - \mathbf{X}_h^{n-1}}{\tau}\cdot (\phi_h\,\boldsymbol\nu_h^n)\, \mathrm{d}\xi = \int_I \partial_{\xi} p_h^n \, \partial_\xi \phi_h \, \rho_h^n \,\mathrm{d}\xi,  && \forall \phi_h\in V_h,  \label{eq:variation_formulation_a}\\
& \int_I \dfrac{\mathbf{X}_h^n - \mathbf{X}_h^{n-1}}{\tau}\cdot (\varphi_h\,\boldsymbol\tau_h^n)\, \mathrm{d}\xi = 0, &&\forall \varphi_h\in V_h,  \label{eq:variation_formulation_b}\\
&  \int_I (p_h^n\,\boldsymbol\nu_h^n + q_h^n\,\boldsymbol\tau_h^n)\cdot \boldsymbol\psi_h\, \mathrm{d}\xi  = \sum_{j = 1}^M\boldsymbol\psi_h(\xi_j)\cdot \int_{t_{n-1}}^{t_{n}}[\boldsymbol\tau_h]_j\, \mathrm{d}t, && \forall \boldsymbol\psi_h\in {\bf V}_h.  \label{eq:variation_formulation_c}
\end{align}
\end{subequations}
The numerical scheme \eqref{eq:variation_formulation} together with \eqref{eq:nu_tau_int}-\eqref{eq:rho_int} is a nonlinear system of equations, which can be solved by using the following Newton's iteration.

\subsection{Properties of the numerical method}
In this subsection, we prove that the proposed numerical scheme \eqref{eq:variation_formulation}, or its equivalent form \eqref{eq:variation_formulation_initial}, can preserve the two geometric structures (i) and (ii) mentioned in the introduction section.


Let $L(t_n)$ be the perimeter of the piecewise linear curve $\Gamma_h(t_n)$, and let $A(t_n)$ be the area of the polygonal region enclosed by $\Gamma_h(t_n)$.

\begin{theorem}\label{MainTHM1}
The solution of \eqref{eq:variation_formulation} has the following properties:
\begin{itemize}
\item[{\rm(1)}]~$A(t_n)= A(t_{n-1})$,
\item[{\rm(2)}]~$L(t_n)\le L(t_{n-1})$ .
\end{itemize}
\end{theorem}

{\it Proof.}$\,$
We use the equivalent formulation \eqref{eq:variation_formulation_initial} to prove these properties.
By integrating the identity
\begin{align*}
\frac{\d }{\d t} A(t)
=-\int_I  \mathbf{v}_h^n \cdot \boldsymbol{\nu}_h \, \d\Gamma_h(t) ,
\end{align*}
in time for $t\in[t_{n-1},t_n]$ and setting $\phi_h=1$ in \eqref{eq:variation_formulation_initial_a}, we obtain
\begin{align*}
A(t_{n}) - A(t_{n-1})
&=-\int_{t_{n-1}}^{t_n} \int_{I} \mathbf{v}_h^n \cdot \boldsymbol{\nu}_h \, \d\Gamma_h(t)\d t \\
&=-\int_{t_{n-1}}^{t_n} \int_{I} \partial_\Gamma p_h^n \cdot \partial_\Gamma 1 \, \d\Gamma_h(t) \d t \\
&=0.
\end{align*}
This proves that the area enclosed by the curve is conserved during the evolution.

We use the following expression for the perimeter of the piecewise linear curve $\Gamma_h(t)$:
$$L(t) = \sum_{j=1}^M \int_{\xi_{j-1}}^{\xi_j} |\partial_\xi \mathbf{X}_h(\xi,t)| \d\xi.$$
Differentiating this expression in time yields
\begin{align*}
\hspace{-10pt}
\frac{\d}{\d t} L(t)
&= \sum_{j=1}^M \int_{\xi_{j-1}}^{\xi_j} \partial_{t} |\partial_\xi \mathbf{X}_h(\xi,t)| \d\xi  \\
&= \sum_{j=1}^M \int_{\xi_{j-1}}^{\xi_j} \frac{\partial_\xi \mathbf{X}_h(\xi,t)}{|\partial_\xi \mathbf{X}_h(\xi,t)|} \cdot\partial_\xi\partial_t \mathbf{X}_h(\xi,t) \d\xi \\
&= \sum_{j=1}^M \int_{\xi_{j-1}}^{\xi_j} {\boldsymbol{\tau}}_h  \cdot\partial_\xi \mathbf{v}_h \, \d\xi \\
&= - \sum_{j=1}^M [{\boldsymbol{\tau}}_h]_j \cdot  \mathbf{v}_h(\xi_j),
\end{align*}
where we have used integration by parts on each subinterval $[\xi_{j-1},\xi_j]$.
By integrating the equality above with respect to time from $t_{n-1}$ to $t_n$, we obtain
\begin{align*}
L(t_n) - L(t_{n-1})
&= - \sum_{j=1}^M \mathbf{v}_h^n(\xi_j) \cdot  \int_{t_{n-1}}^{t_n}[\boldsymbol\tau_h]_j\d t \\
&= - \int_{t_{n-1}}^{t_n} \int_{I} \boldsymbol{\kappa}_h^n \cdot \mathbf{v}_h^n\d\Gamma_h(t) \d t,
\end{align*}
where the last equality is obtained by substituting $\psi_h=\mathbf{v}_h^m$ into \eqref{eq:variation_formulation_initial_c}.
By setting ${\phi}_h = p_h^n$ in \eqref{eq:variation_formulation_initial_a} and $\varphi_h = q_h^n$ in \eqref{eq:variation_formulation_initial_b}, we have
\begin{align*}
&\int_{t_{n-1}}^{t_n}
\int_{I}  \boldsymbol{\kappa}_h^n \cdot \mathbf{v}_h^n
\d\Gamma_h(t) \d t
=
\int_{t_{n-1}}^{t_n}
\int_{I} |\partial_\Gamma (\boldsymbol{\kappa}_h^n\cdot \boldsymbol{\nu}_h)|^2  \d\Gamma_h(t) \d t .
\end{align*}
Finally, by combining the two equalities above, we obtain
\begin{align*}
L(t_n) - L(t_{n-1})
= -\int_{t_{n-1}}^{t_n}
\int_{I} |\partial_\Gamma (\boldsymbol{\kappa}_h^n\cdot \boldsymbol{\nu}_h)|^2  \d\Gamma_h(t) \d t
\le 0 .
\end{align*}
This proves that the perimeter of the computed curve decreases in time.
\qed

\subsection{Newton's iteration}
The nonlinear system \eqref{eq:variation_formulation} can be solved by Newton's iteration. We denote by $(\mathbf{X}_{h,l-1}^n, p_{h,l-1}^n, q_{h,l-1}^n) $ the numerical solution obtained in the $(l-1)^{\rm th}$ iteration, and define
\begin{align}
& \mathbf{X}_{h,l-1}(\xi, t) = \dfrac{t_{n} - t}{\tau}\mathbf{X}_h^{n - 1} + \dfrac{t - t_{n-1}}{\tau}\mathbf{X}_{h,l-1}^n \\
& \boldsymbol\tau_{h,l-1} = \dfrac{\partial_\xi\mathbf{X}_{h,l-1}}{|\partial_\xi\mathbf{X}_{h,l-1}|}
&&\mbox{for $t\in[t_{n-1}, t_{n}]$} .
\end{align}
In the $l^{\rm th}$ iteration, we calculate
\begin{align}
&\boldsymbol\nu_{h,l-1}^n = \int_{t_{n-1}}^{t_{n}} \left(\partial_\xi\mathbf{X}_{h,l-1}\right)^\perp\, \mathrm{d}t \label{eq:newton_nu_tau}\\
&\boldsymbol\tau_{h,l-1}^n = \int_{t_{n-1}}^{t_{n}}\partial_\xi\mathbf{X}_{h,l-1} \, \mathrm{d}t \\
&\rho_{h,l-1}^n = \int_{t_{n-1}}^{t_{n}}\left|\partial_\xi\mathbf{X}_{h,l-1}\right|^{-1}\,\mathrm{d}t \\
& \boldsymbol\theta_{j, l-1}^n = \int_{t_{n-1}}^{t_{n}}\left[\boldsymbol\tau_{h,l-1}\right]_j\,\mathrm{d}t \\
    &\boldsymbol\gamma_{h,l-1}^n = \int_{t_{n-1}}^{t_{n}}\left|\partial_\xi\mathbf{X}_{h,l-1}\right|^{-3}\partial_\xi\mathbf{X}_{h,l-1}\,\mathrm{d}t \\
&a_{h,l-1}^n = \int_{t_{n-1}}^{t_{n}}\left|\partial_\xi\mathbf{X}_{h,l-1}\right|^{-1}\dfrac{t - t_{n-1}}{\tau}\,\mathrm{d}t \\
    &\bm{A}_{h,l-1}^n = \int_{t_{n-1}}^{t_{n}} \left| \partial_\xi\mathbf{X}_{h,l-1}\right|^{-3} \left(\partial_\xi\mathbf{X}_{h,l-1}\right)^{T} \left(\partial_\xi\mathbf{X}_{h,l-1}\right)\,\mathrm{d}t \label{eq:newton_A}
\end{align}
and solve the following linear system to obtain the Newton direction $(\mathbf{X}_\delta, p_\delta, q_\delta)\in {\bf V}_h\times V_h\times V_h$,
\begin{subequations}\label{eq:newton_linear}
\begin{align}
& \int_0^1\dfrac{\mathbf{X}_\delta}{\tau}\cdot\left(\phi_h\boldsymbol\nu_{h,l-1}^n\right)\,\mathrm{d}\xi + \int_0^1\dfrac{\mathbf{X}_{h,l-1}^n -
                 \mathbf{X}_h^{n-1}}{\tau}\cdot\left(\dfrac{\tau}{2}\phi_h(\partial_\xi\mathbf{X}_\delta)^\perp\right)\,\mathrm{d}\xi  \\
&\quad\, - \int_0^1 \partial_\xi p_\delta \partial_\xi\phi_h \rho_{h,l-1}^n\,\mathrm{d}\xi
 + \int_0^1\partial_\xi p_{h,l-1}^n\partial_\xi\phi_h\left(\boldsymbol\gamma_{h,l-1}^n\cdot
                 \partial_\xi\mathbf{X}_\delta\right)\,\mathrm{d}\xi \nonumber \\
& = -\int_0^1\dfrac{\mathbf{X}_{h,l-1}^n - \mathbf{X}_h^{n-1}}{\tau}\cdot \left(\phi_h \boldsymbol\nu_{h,l-1}^n\right) \,
                 \mathrm{d}\xi + \int_0^1\partial_\xi p_{h,l-1}^n \partial_\xi\phi_h\rho_{h,l-1}^n\, \mathrm{d}\xi
                 &&
                 \forall \phi_h\in V_h, \nonumber\\[10pt]
& \int_0^1\dfrac{\mathbf{X}_\delta}{\tau}\cdot\left(\varphi_h\boldsymbol\tau_{h,l-1}^n\right)\,\mathrm{d}\xi +
               \int_0^1\dfrac{\mathbf{X}_{h,l-1}^n - \mathbf{X}_h^{n-1}}{\tau}\cdot\left(\dfrac{\tau}{2}\varphi_h\partial_\xi\mathbf{X}_\delta\right)\,\mathrm{d}\xi  \\
& = - \int_0^1\dfrac{\mathbf{X}_{h,l-1}^n - \mathbf{X}_h^{n - 1}}{\tau}\cdot
               \left(\varphi_h\boldsymbol\tau_{h,l-1}^n\right)\,\mathrm{d}\xi
               &&\forall \varphi_h\in V_h, \nonumber \\[10pt]
& \int_0^1p_\delta\boldsymbol\nu_{h,l-1}^n\cdot\boldsymbol\psi_h\,\mathrm{d}\xi +
               \int_0^1\dfrac{\tau}{2} p_{h,l-1}^n(\partial_\xi\mathbf{X}_\delta)^\perp\cdot\boldsymbol\psi_h\,\mathrm{d}\xi \\
&\quad\,
+ \int_0^1q_\delta\boldsymbol\tau_{h,l-1}^n\cdot\boldsymbol\psi_h \,\mathrm{d}\xi + \int_0^1 \dfrac{\tau}{2} q_{h,l-1}^n\partial_\xi\mathbf{X}_\delta\cdot\boldsymbol\psi_h\,\mathrm{d}\xi \nonumber \\
&\quad\, - \sum_{j = 1}^M\boldsymbol\psi_h(\xi_j)\cdot [\partial_\xi\mathbf{X}_\delta a_{h,l-1}^n]_j + \sum_{j = 1}^M \left[\bm{A}_{h,l-1}^n:\left(\boldsymbol\psi_h(\xi_j)^T\partial_\xi\mathbf{X}_\delta\right)\right]_j \nonumber \\
& = - \int_0^1\left(p_{h,l-1}^n\boldsymbol\nu_{h,l-1}^n + q_{h,l-1}^n
                \boldsymbol\tau_{h,l-1}^n\right)\cdot\boldsymbol\psi_h\,\mathrm{d}\xi + \sum_{j = 1}^{M}\boldsymbol\psi_h(\xi_j)\cdot\boldsymbol\theta_{j, l-1}^n && \forall \boldsymbol\psi_h\in {\bf V}_h.
\nonumber
\end{align}
\end{subequations}
The above linearized problem \eqref{eq:newton_linear} is obtained by using the first-order Taylor expansion at the point $(\mathbf{X}_{h,l-1}^n, p_{h,l-1}^n, q_{h,l-1}^n)$ of the nonlinear system \eqref{eq:variation_formulation}, and by denoting $\mathbf{X}_{\delta} = \mathbf{X}^n - \mathbf{X}_{h,l-1}^n$, $p_\delta = p^n - p_{h,l-1}^n$ and $q_\delta = q^n - q_{h,l-1}^n$.

The algorithm for solving \eqref{eq:variation_formulation} by Newton's iteration is stated as follows:
\begin{itemize}
\item Step 1. Choose three proper tolerances ${\rm TOL}, {\rm TOL'}, {\rm TOL''}>0$. Take the initial guess of Newton's iteration to be
$$
\mathbf{X}_{h,0}^n = \mathbf{X}_h^{n-1},\quad
p_{h,0}^n = p_h^{n-1}
\quad\mbox{and}\quad
q_{h,0}^n = q_h^{n-1} .
$$

\item Step 2.
For given $(\mathbf{X}_{h,l-1}^{n-1},p_{h,l-1}^{n-1},q_{h,l-1}^{n-1})$, calculate the quantities in \eqref{eq:newton_nu_tau}-\eqref{eq:newton_A} and solve the linear system \eqref{eq:newton_linear} to obtain $(\mathbf{X}_\delta, p_{\delta}, q_{\delta})\in {\bf V}_h\times V_h\times V_h$. Then set
$$
\mathbf{X}_{h,l}^n = \mathbf{X}_{h,l-1}^n + \mathbf{X}_\delta,
\quad
p_{h,l}^n = p_{h,l-1}^n + p_\delta
\quad\mbox{and}\quad
q_{h,l}^n = q_{h,l-1}^n + q_\delta .
$$

\item Step 3. If $\|\mathbf{X}_{h,l}^n - \mathbf{X}_{h,l-1}^n\|_{L^\infty}\le {\rm TOL}$, $\|p_{h,l}^n - p_{h,l-1}^n\|_{L^\infty} \le {\rm TOL}'$ and $\|q_{h,l}^n - q_{h,l-1}^n\|_{L^\infty} \le {\rm TOL}''$, then set
$$
\mathbf{X}_h^n = \mathbf{X}_{h,l}^n,
\quad
p_h^n = p_{h,l}^n
\quad\mbox{and}\quad
q_h^n = q_{h,l}^n
$$
and stop; otherwise, set $l = l + 1$ and goto Step 2.

\end{itemize}

\section{Numerical results}

In this section, we present several numerical experiments to demonstrate the accuracy of the proposed numerical scheme \eqref{eq:variation_formulation} and the effectiveness of the method in preserving the two geometric structures.

\subsection{Convergence test}

We consider the surface diffusion flow of a closed curve with initial shape being an ellipse
$$x^2+4y^2=4.$$
We approximate the surface diffusion flow by the proposed numerical scheme \eqref{eq:variation_formulation} and present the relative errors in Table \ref{tab:newton}, with
\begin{equation}
    e_{h, \tau}(t) = \|\mathbf{X}_{h,\tau} - \mathbf{X}_{h/2, \tau/4}\|_{L^{\infty}(I)} ,
\end{equation}
where $\mathbf{X}_{h, \tau}$ is the numerical solution obtained with mesh size $h:=1/M$ and time stepsize $\tau$. The order of convergence is calculated by
\begin{equation*}
\mbox{order of convergence}
= \log \bigg(\frac{\|\mathbf{X}_{h,\tau} - \mathbf{X}_{h/2, \tau/4}\|_{L^{\infty}(I)}}{\| \mathbf{X}_{h/2, \tau/4} - \mathbf{X}_{h/4, \tau/16}\|_{L^{\infty}(I)}}\bigg)\bigg/\log(2) ,
\end{equation*}
based on the finest three meshes. The nonlinear system is solved by using Newton's iteration with a tolerance error of $10^{-10}$.

From Table \ref{tab:newton} we can see that the proposed numerical scheme has second-order convergence in space when $\tau=O(h^2)$. This is the same as the semi-implicit parametric FEM considered in \cite{Bao17} (see \cite[Table 1 on p. 389]{Bao17}).

\begin{table}[hbt!]
\begin{center}
\caption{Convergence rates in the $L^\infty$ norm, with $h_0 = 1/8$ and $\tau_0 = 0.04$. The initial curve is an ellipse $x^2+4y^2=4$.}\label{tab:newton}
\vspace{0cm}{\setlength{\tabcolsep}{7pt}
	\centering
	\begin{tabular}{ |l|c|c|c|c|c|cc|c|}
		\hline
	  & $\begin{array}{l}
	   h=1/8 \\
	   \tau=1/25
	   \end{array}$ & $\begin{array}{l}
	   h=1/16 \\
	   \tau=1/100
	   \end{array}$ & $\begin{array}{l}
	   h=1/32 \\
	   \tau=1/400
	   \end{array}$ & $\begin{array}{l}
	   h=1/64 \\
	   \tau=1/1600
	   \end{array}  $ & order \\
		\hline
$t_N=0.2$    & 1.10E-2 & 3.69E-3 & 9.96E-4  & 2.55E-4  & 1.97  \\
    \hline
$t_N=0.5$   & 2.42E-2 & 7.15E-3 & 1.93E-3 & 5.07E-4  & 1.93  \\
    \hline
$t_N=2$   & 1.41E-2 & 4.08E-3 & 1.03E-3 & 2.57E-4  & 2.00  \\
    \hline
   \end{tabular}}
\end{center}
\end{table}

\begin{figure}[htp!]
    \centering
    \includegraphics[width = .5\textwidth]{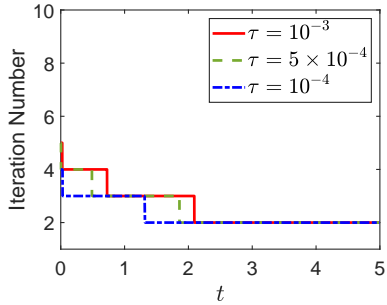}
    \vspace{-5pt}
    \caption{The number of Newton's iterations in each time step.} \label{fig:iteration_number}
\end{figure}

Figure \ref{fig:iteration_number} shows the number of Newton's iterations in the evolution process, with three different time stepsizes. One can see that the method needs only a few iterations to attain the desired accuracy, and the time stepsize does not have much influence on the numerical of iterations.

\subsection{Structure preservation}
%

\begin{figure}[htbp!]
    \centering
    \includegraphics[width = .95\textwidth]{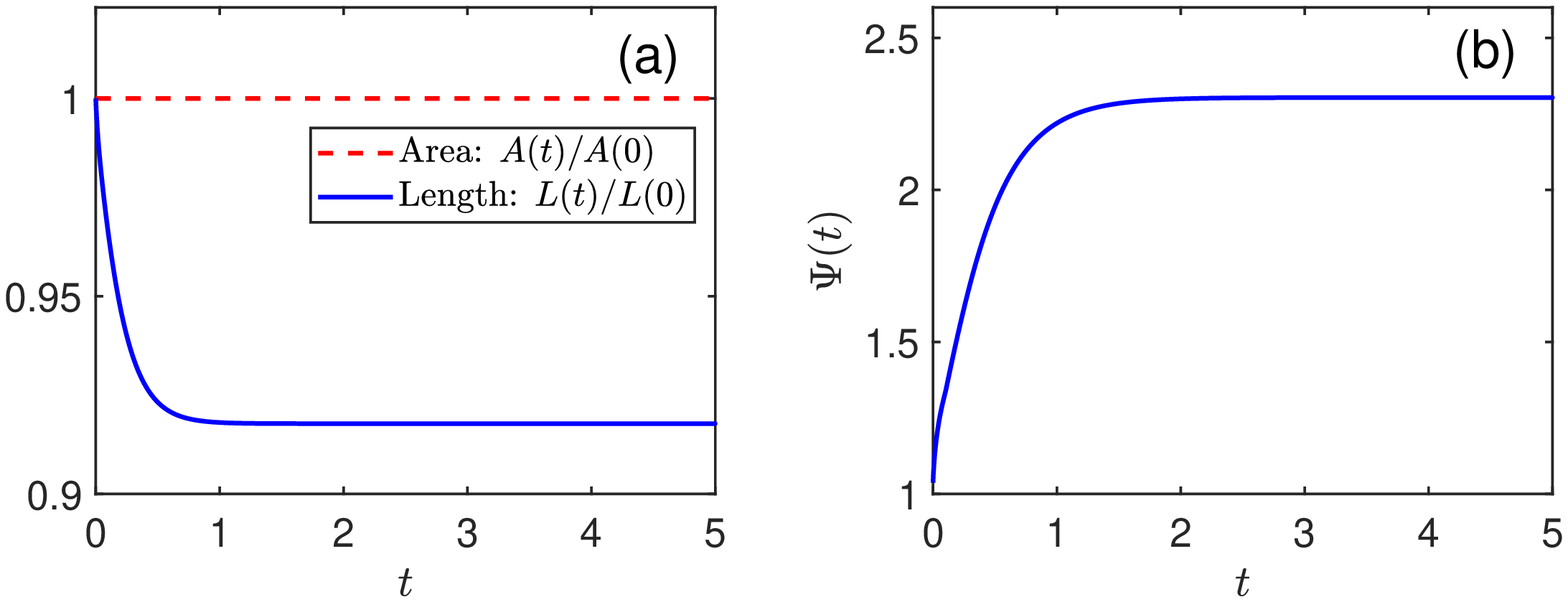}
    \vspace{-5mm}
  \caption{Numerical solution given by the proposed method with $h = 1/32$ and $\tau = 10^{-4}$. \newline
\indent\hspace{45pt}  (a) The normalized perimeter and normalized area enclosed by the curve; \newline
\indent\hspace{45pt}  (b) The mesh distribution function $\Psi(t)$.}
  \label{fig:numerical_property_newton}
\end{figure}

For the example considered in the last subsection, we present in
Figure \ref{fig:numerical_property_newton} the evolution of the normalized perimeter $L(t_{n})/L(0)$, the normalized enclosed area $A(t_{n})/A(0)$, and a mesh distribution function $\Psi(t_{n})$, defined by~\cite{Bao17}
\begin{align*}
\Psi(t=t_{n}) = \dfrac{\max\limits_{1\le j\le M}\left|\mathbf{x}_j^n - \mathbf{x}_{j - 1}^n\right|}{\min\limits_{1\le j\le M}\left|\mathbf{x}_j^n - \mathbf{x}_{j-1}^n\right|} .
\end{align*}

\begin{figure}[htp!]
\centering
\includegraphics[width = 0.95\textwidth]{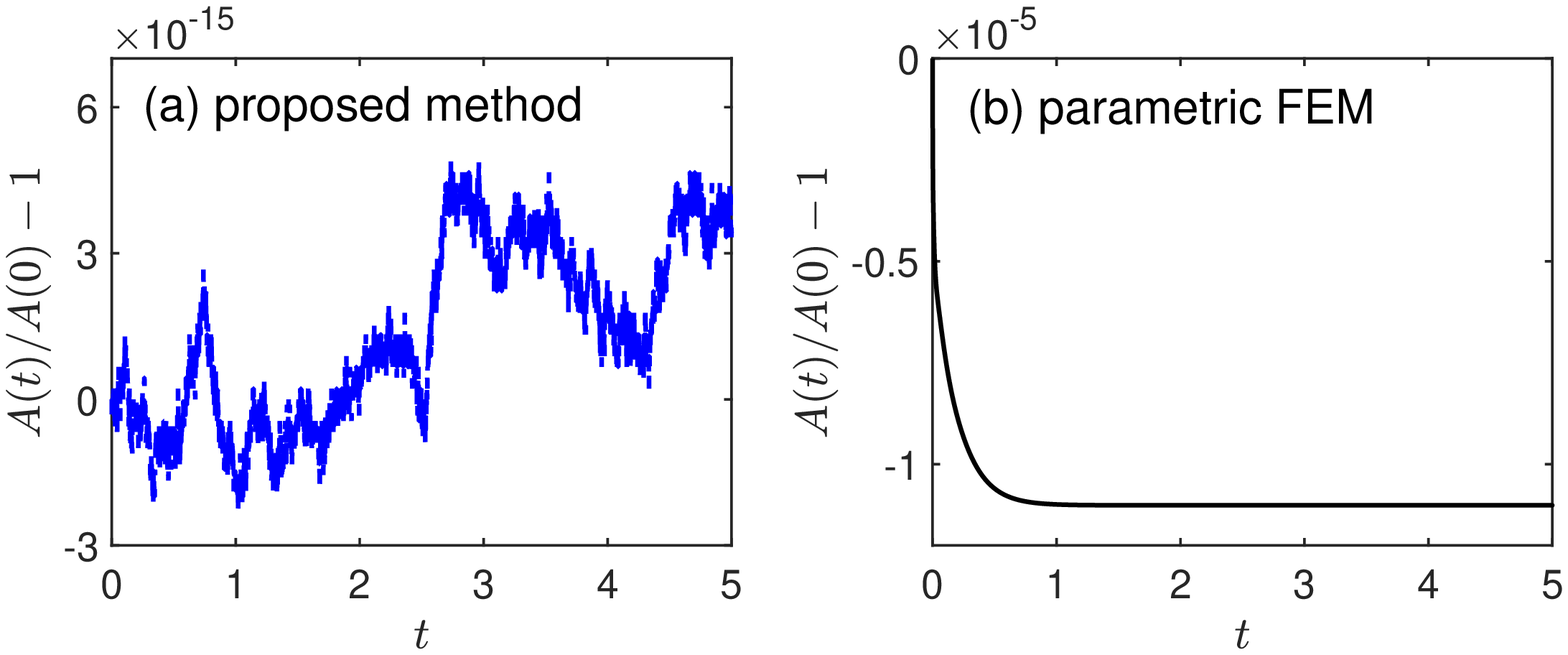}
\vspace{-5mm}
\caption{Comparison between the proposed method and the parametric FEM in \cite{Bao17}. \newline
\indent\hspace{45pt}  Both methods use $h = 1/32$ and $\tau = 10^{-4}$.} \label{fig:area_accurate}
\end{figure}

From Figure \ref{fig:numerical_property_newton}, one can see that the proposed method can preserve the conservation of $A(t)/A(0)$ and the decrease of $L(t)/L(0)$, while the mesh distribution function $\Psi(t)$ grows from $1$ to around $2.3$.
Figure \ref{fig:area_accurate} (a) shows that the area enclosed by the curve is conserved with machine precision by the proposed method; Figure \ref{fig:area_accurate} (b) shows that the parametric FEM considered in \cite{Bao17} and Barrett et al. \cite{Barrett07a} does not strictly preserve the area conservation.

\subsection{Evolution of a ``flower'' shape}

\begin{figure}[htp!]
\centering
\includegraphics[width = .95\textwidth]{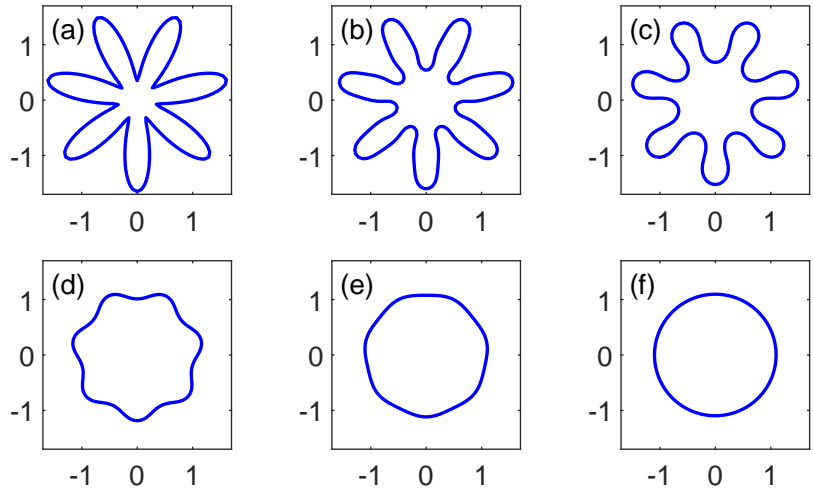}
\caption{Several steps in the evolution of an initially ``flower'' shape toward its equilibrium at different times: (a) $t = 0$; (b) $t = 10^{-4}$; (c) $t = 0.001$; (d) $t = 0.005$; (e) $t = 0.006$ and (f) $t = 0.01$, where $M = 210$ and $\tau = 10^{-6}$.} \label{fig:evolution_7flower}
\end{figure}

\begin{figure}[htp!]
    \centering
    \includegraphics[width = .90\textwidth]{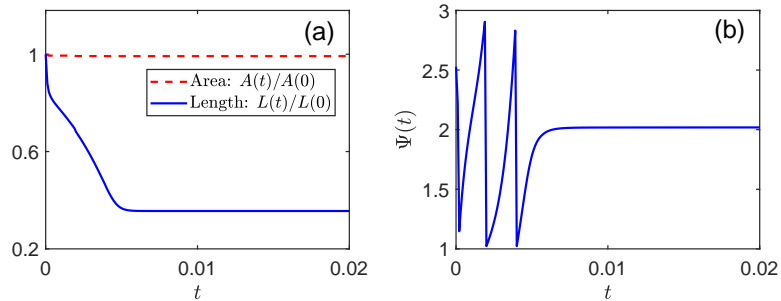}
  \caption{The corresponding temporal evolutions shown in Fig.~\ref{fig:evolution_7flower} for: (a) the normalized total free energy and the normalized area; (b) the mesh distribution function $\Psi(t)$, respectively.}
  \label{fig:numerical_property_7flower}
\end{figure}

We apply our method to the numerical simulation of surface diffusion flow which initially is a ``flower'' shape, given in the polar coordinate by
\begin{equation}
\left\{
\begin{aligned}
&x = \left[1 + 0.65\sin(7\theta)\right]\cos\theta,\\
&y = \left[1 + 0.65\sin(7\theta)\right]\sin\theta,
\end{aligned}
\right. \qquad  \theta\in[0, 2\pi].
\end{equation}
The shapes of the curves are presented in Figure \ref{fig:evolution_7flower} at several different time levels.
The evolution of the normalized area, normalized perimeter, and mesh distribution function is presented in Figure \ref{fig:numerical_property_7flower}.
The two figures show that the flower-shape curve evolves gradually to a circle with the same area and shorter perimeter.


\section{Conclusion}

We have proposed a novel numerical method for simulating the isotropic surface diffusion flow of closed curves in two dimensions. Based on a new weak formulation of the surface diffusion equation, we have introduced a piecewise linear finite element discretization for the weak formulation to obtain a nonlinear system of equations, which can be solved by Newton's iteration.
We have rigorously proved that the proposed numerical method can simultaneously preserve the two geometric structures, i.e., area conservation and perimeter decrease. Through numerical experiments it is shown that the proposed method is second-order convergent in the $L^{\infty}$-norm, and can retain the area conservation and perimeter decrease with machine precision.
For a complex flower-shape curve,
we have obtained satisfactory numerical results by combining the proposed method with a mesh redistribution technique.

\bibliographystyle{plain}
\bibliography{mybib}

\end{document}